\documentclass[a4paper,english, nostix]{birthday}
\usepackage{amsmath,babel}

\newcommand{\R}{\mathbb{R}}

\newcommand{\DD}{\mathcal{D}}
\newcommand{\Z}{\mathbb{Z}}
\newcommand{\N}{\mathbb{N}}
\newcommand{\C}{\mathbb{C}}
\newcommand{\Q}{\mathbb{Q}}
\newcommand{\LL}{\mathcal{L}}

\newcommand{\E}{{\mathcal E}}

\newcommand{\po}{\partial}

\newcommand{\ve}{\varepsilon}

\newcommand{\la}{\langle}
\newcommand{\ra}{\rangle}

\newcommand{\loc}{{\text{\rm loc}}}
\newcommand{\X}{\times}

\renewcommand{\d}{\delta}
\renewcommand{\l}{\lambda}

\renewcommand{\a}{\alpha}
\renewcommand{\b}{\beta}

\newcommand{\s}{\sigma}
 
\newcommand{\z}{\zeta}
\renewcommand{\k}{\kappa}
\newcommand{\sgn}{\text{\rm sgn}}

\newcommand{\Om}{\Omega}

\newcommand{\supp}{\text{\rm supp}\,}

\newcommand{\M}{{\mathcal M}}

\renewcommand{\E}{{\mathcal E}}
\newcommand{\DM}{\mathcal D\mathcal M}

\renewcommand{\supp}{\text{\rm supp}\,}

\newcommand{\bff}{{\mathbf f}}
\newcommand{\bfa}{{\mathbf a}}
\newcommand{\bfg}{{\mathbf g}}

\newcommand{\esslim}{\operatorname{ess}\!\lim}

\newcommand{\AP}{\operatorname{AP}}
\newcommand{\BAP}{\operatorname{BAP}} 
\newcommand{\SAP}{\operatorname{SAP}} 
\newcommand{\Me}{\operatorname{M}}
\newcommand{\Sp}{\operatorname{Sp}}
\newcommand{\Gr}{\operatorname{Gr}}
\renewcommand{\subset}{\subseteq}

\newcommand{\ut}{\underline{t}}
\newcommand{\ux}{\underline{x}}

\newcommand{\intl}{\int\limits}
\newcommand{\iintl}{\iint\limits}

\newtheorem{definition}{Definiton} 
\newtheorem{theorem}{Theorem}
\newtheorem{lemma}{Lemma}
\newtheorem{proposition}{Proposition} 

\newcommand{\medint}{{\mbox{\vrule height3.5pt depth-2.8pt
          width4pt}\mkern-13mu\int\nolimits}}
\newcommand{\Medint}{\mkern12mu\mbox{\vrule height4pt
         depth-3.2pt
          width5pt}\mkern-16.5mu\int\nolimits}

\title{On the Decay of Almost Periodic Solutions for Certain \\ Degenerate Parabolic-Hyperbolic Equations}
\author{Hermano Frid \\ \\ \\ \rightline{\sl Dedicated to Helge Holder on his 60th birthday.}}
\contact[hermano@impa.br]{Instituto de Matem\' atica Pura e Aplicada-IMPA\\ Estrada Dona Castorina, 110, CEP 22460-320, Rio de Janeiro, RJ, Brazil}

\begin{document}
\begin{abstract}
We discuss the  well-posedness and decay of Besicovitch almost periodic solutions for a class of nonlinear degenerate anisotropic hyperbolic-parabolic equations. 
In our definition of weak entropy solution the initial data is only assumed in a weak sense, and in this connection the strong trace of the solution  in the initial time hyperplane is also established.  

\end{abstract}


\section{Introduction}\label{S:1}

We consider the following Cauchy problem for a nonlinear degenerate anisotropic hyperbolic-parabolic equation
\begin{align}
&\po_t u+\nabla\cdot \bff(u) =\Delta_{x''} b(u), \quad \text{in $(0,\infty)\X\R^d$}, \label{e1.1} \\
&u(0,x)=u_0(x),\label{e1.2}
\end{align}
where $\bff:\R\to\R^d$, $\bff=(f_1,\cdots,f_d)$,  and $b:\R\to\R$ are smooth functions,  $b'(u)\ge0$,  we split $\R^d=\R^{d'}\X\R^{d''}$, for $x\in\R^d$, we write $x=(x',x'')$, $x'\in\R^{d'}$, $x''\in\R^{d''}$, and $\Delta_{x''}=\po_{x_{d'+1}}^2+\cdots+\po_{x_{d}}^2$.

We assume the following non-degeneracy condition. Denoting $\bfa(\xi)=\bff'(\xi)$, $\bfa(\xi)=(a_1(\xi),\cdots,a_d(\xi))$, $\pi_{d'}(\bfa)=(a_1,\cdots,a_{d'})$,  we assume that, for any $(\z_0,\z)\in\R^{d'+1}$, $(\z_0,\z)\ne(0,0)$, 
 \begin{equation}\label{e1.3}
 \LL^1\bigl\{\xi\in\R\,:\, \z_0+ \pi_{d'}(\bfa)(\xi)\cdot\z=0\}=0,
 \end{equation} 
 where $\LL^1$ is the Lebesgue measure on $\R$, and we also assume that
 \begin{equation}\label{e1.4}
 \LL^1\{\xi\in\R\,:\, b'(\xi)=0\}=0.
 \end{equation}

In this paper, we are interested in analyzing the asymptotic behavior of solutions of \eqref{e1.1},\eqref{e1.2} with initial data satisfying  
\begin{equation}\label{e1.u0}
u_0\in L^\infty(\R^d)\cap \BAP(\R^d).
\end{equation}
Here, $\BAP(\R^d)$ denotes the space of the  Besicovitch almost periodic functions (with exponent $p=1$), which can be defined as the completion of the space of trigonometric polynomials, i.e., finite sums $\sum\limits_{\l} a_{\l} e^{2\pi i\l\cdot x}$ ($i=\sqrt{-1}$ is the purely imaginary unity)   under the semi-norm
$$
N_1(g)=\limsup_{R\to\infty}\frac1{R^d}\int_{C_R}|g(x)|\,dx,
$$
where, for $R>0$,
$$
C_R:=\{x\in\R^d\,:\, |x|_\infty:=\max_{i=1,\cdots,d}|x_i|\le R/2\}.
$$
We observe that the semi-norm $N_1$ is indeed a norm over the trigonometric polynomials, so the referred completion through it is well defined Banach space. Equivalently, the space $\BAP(\R^d)$ is also the completion through $N_1$ of the space of uniform (or Bohr) almost periodic functions, $\AP(\R^d)$, which the defined as the closure in the $\sup$-norm of the trigonometric polynomials. 

We begin with the definition of entropy solution of  the Cauchy problem \eqref{e1.1},\eqref{e1.2}.  We use the notation $\nabla_{x''}:=(\underset{d'}{\underbrace{0,\cdots,0}},\po_{x_{d'+1}},\cdots,\po_{x_d})$, and also $\nabla_{x'}:=\nabla_x-\nabla_{x''}$. 

\begin{definition} \label{D:1.1} A function $u\in L^\infty((0,\infty)\X\R^d)$ is a weak entropy solution to  the problem \eqref{e1.1},\eqref{e1.2} in $\R^{d+1}_+:=(0,\infty)\X\R^d$, if it satisfies:
\begin{enumerate}
\item[(i)] (Regularity) For any $R>0$ and $T>0$,  if $C_{R,T}:= (0,T)\X C_R$, we have
\begin{equation}\label{e1.5}
\nabla_{x''} b(u) \in L^2(C_{R,T}),
\end{equation}
and, for some constant $C(T)>0$, we have
\begin{equation}\label{e1.5'}
\limsup_{R\to+\infty} R^{-d}\iintl_{C_{R,T}} |\nabla_{x''}b(u)|^2\,dx\,dt \le C(T).
\end{equation}

\item[(i)] (Entropy condition) For any $\phi\in C_0^\infty(\R_+^{d+1})$, $\phi\ge0$, and any $k\in\R$, we have
\begin{equation}\label{e1.6}
\iintl_{\R^{d+1}_+}|u-k|\phi_t+\sgn(u-k)(\bff(u)-\bff(k))\cdot\nabla\phi -\nabla_{x''}| b(u)-b(k)|\cdot\nabla_{x''}\phi\,dx\,dt \ge0.
\end{equation}

\item[(ii)] (Initial condition) For any $\varphi\in C_0(\R^{d+1})$ we have
\begin{multline}\label{e1.7}
\iintl_{\R_+^{d+1}} u\varphi_t + \bff(u)\cdot\nabla\varphi -\nabla_{x''}b(u)\cdot \nabla_{x''}\varphi \, dx\,dt\\
+\intl_{\R^d}u_0(x)\varphi(0,x)\,dx=0.
\end{multline}
\end{enumerate} 
\end{definition}  

The first main result of this paper is the following.

\begin{theorem}\label{T:1.1} Let $u$ be weak entropy solution of \eqref{e1.1},\eqref{e1.2}. Then, 
$$
u\in C([0,\infty), L_{\loc}^1(\R^d)). 
$$
In particular, for any $R>0$, 
\begin{equation}\label{e1.7'}
\lim_{t\to0+} \int_{|x|<R} |u(t,x)-u_0(x)|\,dx=0.
\end{equation}

\end{theorem}

For any $g\in\BAP(\R^d)$, its mean value $\Me(g)$, defined by
$$
\Me(g):=\lim_{R\to\infty} R^{-d}\int_{C_R} g(x)\,dx,
$$
exists (see, e.g., \cite{B}). The mean value $\Me(g)$ is also denoted by $\medint_{\R^d}g\,dx$.  Also, the  Bohr-Fourier coefficients of $g\in\BAP(\R^d)$ 
$$
a_\l= \Me(g e^{-2\pi i\l\cdot x}),
$$ 
are well defined and we have that the spectrum of $g$, defined by
$$
\Sp(g):=\{ \l\in\R^g \,:\, a_{\l}\ne 0\},
$$
is at most countable (see, e.g., \cite{B}). We denote by  $\Gr(g)$ the smallest additive subgroup of $\R^d$ containing $\Sp(g)$ ({\em cf.}\ \cite{Pv2}, where $\Gr(g)$ was introduced and  denoted by $M(g)$). 

The second main result of this paper is the following.

\begin{theorem}\label{T:1.2}  For any $u_0\in L^\infty(\R^d)$, there exists a unique weak entropy solution $u(t,x)$ of \eqref{e1.1},\eqref{e1.2}. Moreover, if $u_0$ satisfies \eqref{e1.u0}, then 
\begin{equation}\label{e1.8}
u \in  C((0,\infty),\BAP(\R^d))\bigcap L^\infty(\R_+^{d+1}), 
\end{equation}
and $\Gr(u(t,\cdot))\subset \Gr(u_0)$. Further, 
\begin{equation}\label{e1.8'}
\lim_{t\to+\infty}\Me(|u(t,\cdot)-\Me(u_0)|)=0.
\end{equation}
\end{theorem}

The proofs of Theorem~\ref{T:1.1} and Theorem~\ref{T:1.2} are given in the following sections. 

As a final statement in this introduction, we want to establish a result which is the analog for \eqref{e1.1} of proposition~1.3 in \cite{Pv2}, for scalar conservation laws.  We will need the following  technical lemma of  \cite{Pv2}, to which we refer for the proof.

\begin{lemma} \label{L:1.1} Suppose that $u(x,y)\in L^\infty(\R^n\X\R^m)$, 
$$
E=\{x\in\R^n\,:\, \text{ $(x,y)$ is a Lebesgue point of $u(x,y)$ for a.e.\ $y\in\R^m$} \}.
$$
Then $E$ is a set of full measure and $x\in E$ is a common Lebesgue point of the functions $I(x)=\intl_{\R^m}u(x,y)\rho(y)\,dy$, for all $\rho\in L^1(\R^m)$.
\end{lemma}

\begin{proposition} \label{P:1.1} {\bf(mean $L^1$-contraction property).} Let $u(t,x),v(t,x) \in L^\infty(\R_+^{d+1})$ be two weak entropy solutions of \eqref{e1.1},\eqref{e1.2}, with initial data $u_0,v_0\in L^\infty(\R^d)$. Then for a.e.\ $0< t_0<t_1$ 
\begin{equation}\label{e1.9}
N_1(u(t_1,\cdot)-v(t_1,\cdot))\le N_1(u(t_0,\cdot)-v(t_0,\cdot)),
\end{equation}
and also for a.e.\ $t>0$, 
\begin{equation}\label{e1.9'}
N_1(u(t,\cdot)-v(t,\cdot))\le N_1(u_0-v_0),
\end{equation}
\end{proposition}

\noindent
\textsc{Proof:}  We follow closely  with the due adaptations the proof of proposition~1.3 in \cite{Pv2}. We first recall that by using the doubling of variables method of Kruzhkov \cite{Kr}, as adapted by Carrillo \cite{Ca} to the isotropic degenerate parabolic case, we obtain
\begin{equation}\label{e1.10}
|u-v|_t+\nabla\cdot \sgn(u-v)(\bff(u)-\bff(v))\le \Delta_{x''}|b(u)-b(v)|,
\end{equation}
in the sense of distributions in $\R_+^{d+1}$. As usual, we define a sequence approximating the indicator function of the interval $(t_0,t_1]$ , by setting for $\nu\in\N$,
$$
\d_\nu(s)=\nu\rho(\nu s),\quad \theta_\nu(t)=\int_0^t\d_\nu(s)\,ds=\int_0^{\nu t}\rho(s)\,ds,
$$
where $\rho\in C_0^\infty(\R)$, $\supp \rho\subset [0,1]$, $\rho\ge0$, $\int_{\R}\rho(s)\,ds=1$.   We see that $\d_\nu(s)$ converges to the Dirac me sure in the sense of distributions in 
$\R$ while $\theta_\nu(t)$ converges everywhere to the Heaviside function. For $t_1>t_0>0$, if $\chi_\nu(t)=\theta_\nu(t-t_0)-\theta_\nu(t-t_1)$, then $\chi_\nu\in C_0^\infty(\R_+)$,
$0\le\chi_\nu\le 1$, and the sequence $\chi_\nu(t)$ converges everywhere, as $\nu\to\infty$, to the indicator function of the interval $(t_0,t_1]$. Let us take $g\in C_0^\infty(\R^d)$, satisfying $0\le g\le 1$, $g(y)\equiv 1$ in the cube $C_1$, $g(y)\equiv 0$ outside the cube $C_k$, with $k>1$. We apply \eqref{e1.10} to the yest function $\varphi=R^{-d}\chi_\nu(t)g(x/R)$, for $R>0$. We then get
\begin{multline}\label{e1.11}
\int_0^\infty\bigl(R^{-d}\intl_{\R^d}|u(t,x)-v(t,x)|g(x/R)\,dx\bigr)(\d_\nu(t-t_0)-\d_\nu(t-t_1))\,dt\\
+R^{-d-1}\iint_{\R_+^{d+1}}\sgn(u-v)(\bff(u)-\bff(v))\cdot\nabla_yg(x/R)\chi_\nu(t)\,dx\,dt\\
-R^{-d-1}\iint_{R_+^{d+1}}\nabla_{x''}|b(u)-b(v)|\cdot\nabla_{y''}g(x/R)\chi_\nu(t)\,dx\,dt\ge0.
\end{multline}
Define
$$
F=\{t>0\,:\, \text{$(t,x)$ is a Lebesgue point of $|u(t,x)-v(t,x)|$ for a.e.\ $x\in\R^d$} \}.
$$
As a consequence of Fubini's theorem, $F$ is a set of full Lebesgue measure and by Lemma~\ref{L:1.1} each $t\in F$ is a Lebesgue point of the functions
$$
I_R(t)=R^{-d}\int_{\R^d}|u(t,x)-v(t,x)| g(x/R)\,dx,
$$
for all $R>0$ and all $g\in C_0(\R)$. Now we assume $t_0,t_1\in F$ and take the limit as $\nu\to\infty$ in \eqref{e1.11} , to get
\begin{multline}\label{e1.12}
I_R(t_1)\le I_R(t_0)+ R^{-d-1}\iintl_{(t_0,t_1)\X\R^d}\sgn(u-v)(\bff(u)-\bff(v))\cdot\nabla_y g(x/R)\,dx\,dt\\
-R^{-d-1}\iintl_{(t_0,t_1)\X\R^d} \nabla_{x''}|b(u)-b(v)|\cdot\nabla_{y''} g(x/R)\,dx\,dt.
\end{multline}
Now, we have
\begin{multline}\label{e1.13}
R^{-d-1}\bigl |\iintl_{(t_0,t_1)\X\R^d}\sgn(u-v)(\bff(u)-\bff(v))\cdot \nabla_y g(x/R)\,dx\,dt\bigr|\\
\le R^{-1} \|\bff(u)-\bff(v)\|_\infty\iintl_{(t_0,t_1)\X\R^d}|\nabla_y g(y)|\,dy\,dt\to 0,\quad \text{as $R\to\infty$}.
\end{multline}
Also, we have
\begin{multline}\label{e1.14}
R^{-d-1}\bigl|\iintl_{(t_0,t_1)\X\R^d}\nabla_{x''}|b(u)-b(v)|\cdot\nabla_{y''} g(x/R)\,dx\,dt\bigr|\\
\le R^{-1}(R^{-d}\iintl_{(t_0,t_1)\X C_{kR}}|\nabla_{x''}|b(u)-b(v)||^2\,dx\,dt)^{1/2}(\iintl_{(t_0,t_1)\X\R^d}|\nabla_y g(y)|^2\,dy\,dt)^{1/2}\\
\longrightarrow0\qquad  \text{as $R\to\infty$},
\end{multline}
 where we have used  \eqref{e1.5'}. On the other hand, we have
 $$
 N_1(u(t,\cdot)-v(t,\cdot))\le \limsup_{R\to\infty} I_R(t)\le k^d N_1(u(t,\cdot)-v(t,\cdot)),
 $$
 so taking the limit as $R\to\infty$ in \eqref{e1.12}, for $t_0, t_1\in F$, $t_0<t_1$,   we get 
 $$
 N_1(u(t_1,\cdot)-v(t_1,\cdot))\le k^{d} N_1(u(t_0,\cdot)-v(t_0,\cdot)),
 $$
 and since $k>1$ is arbitrary we can make $k\to 1+$ to get the desired result. Finally, for $t_0=0$, we use Theorem~\ref{T:1.1} to send $t_0\to0+$ in \eqref{e1.12} and proceed exactly as we have just done. 
  
  $\hfill\square$
  
There is a large literature related with degenerate parabolic equations, being the first important contribution by  Vol'pert and Hudjaev in \cite{VH}. Uniqueness for the homogeneous Dirichlet problem, for the isotropic case, was only achieved  many years later by Carrillo in \cite{Ca}, using an extension of Kruzhkov's doubling of variables method \cite{Kr}. The result in \cite{Ca} was  extended to non-homogeneous Dirichlet data by Mascia, Porretta and Terracina in \cite{MPT}.   An $L^1$ theory for the Cauchy problem for anisotropic degenerate parabolic equations was established by Chen and Perthame \cite{CP}, based on the kinetic formulation (see \cite{PB}), and later also obtained using Kruzhkov's approach in \cite{CK} (see also, \cite{KR}, \cite{EJ} and the references therein). Decay of almost periodic solutions for general nonlinear systems of conservation laws of parabolic and hyperbolic types was first addressed in \cite{Fr0}, as an extension of the ideas put forth in \cite{CF0}. Only recently the problem of the decay of almost periodic solutions was retaken, specifically for scalar conservation laws, by Panov in \cite{Pv2}, where some elegant ideas were introduced to successfully extend the result in \cite{Fr0} in that specific case. Here we use the elegant ideas in \cite{Pv2}, combined either with \cite{CP2} or \cite{Fr0}, to prove our decay result.  
   
We also remark that all the results in this paper hold with almost no change in the proofs if instead of $\nabla_{x''} b(u)$ we have $\po^2_{x_{d'+1}}b^{d'+1}(u)+\cdots+\po_{x_d}^2b^d(u)$, for $b^{d'+1}(u),\cdots,b^{d}(u)$ satisfying the same conditions as $b(u)$. 
   
The remaining discussion in this paper is organized as follows. In Section~\ref{S:2}, the proof Theorem~\ref{T:1.1} is given. Section~\ref{S:3} contains the proof of Theorem~\ref{T:1.2}. 

\section{Strong trace at hyperplanes  $t=t_0\ge0$.}\label{S:2}

In this section we prove Theorem~\ref{T:1.1}. This amounts to proving the strong trace property for the weak entropy solution of \eqref{e1.1},\eqref{e1.2},  at all hyperplane $t=t_0$, for all $t_0\ge0$. Indeed, by the Gauss-Green Theorem (see, e.g., \cite{CF2}, \cite{CF3}), applied to the (divergence-free) $L^2$-divergence-measure field $(u, \bff(u)-\nabla_{x''}b(u))$, we easily deduce that the limits $\lim_{t\to t_{0\pm}} u(t,x)$ exist in  the weak star topology of $L^\infty(\R^d)$, for $t_0>0$, and just the limit for $t_{0+}$ when $t_0=0$. By the same result, for $t_0>0$, using the fact that the referred field is divergence-free, we easily deduce that the limits for $t_{0+}$ and $t_{0-}$ must coincide.  

We recall that for $Q=(0,T)\X\Om'\X\Om''$, $\Om'\subset\R^{d'}$, $\Om''\subset\R^{d''}$, bounded open sets with smooth boundaries, in $\R^{d'}$ and $\R^{d''}$, and the lateral boundary $\Gamma':=(0,T)\X\po\Om'\X\Om''$, the strong trace property at $\Gamma'$ for any function $u\in L^\infty(Q)$ satisfying \eqref{e1.5} and \eqref{e1.6}, for any $0\le\phi\in C_0^\infty(Q)$, has been proved in \cite{FL}. The proof in \cite{FL} can be easily adapted to give the proof of the strong trace property at the hyperplanes $t=t_0$. However, the latter is actually simpler and we outline a direct proof here for the reader's convenience. 

We start by recalling the kinetic formulation  (see, e.g., \cite{CP}) that if $u\in L^\infty(\R_+^{d+1})$ satisfies \eqref{e1.5} and \eqref{e1.6}, then the function
 $$
 f(t,x,\xi)=\chi(\xi; u(x,t)),\quad \text{where}\quad \chi_u(\xi;u):=\begin{cases} -1, \quad u\le \xi< 0,\\ 1,\quad  0<\xi \le u,\\ 0,\quad |\xi|>|u|,\end{cases}
 $$
 satisfies
 \begin{equation} \label{e2.1}
 \po_t f+\bfa(\xi)\cdot\nabla f-b'(\xi)\Delta_{x''} f=\po_{\xi} m,
 \end{equation}
 in the sense of distributions in $\DD'(\R_+^{d+1}\X(-L,L))$, with $\bfa(\xi)=\bff'(\xi)$, for some $m\in\M_\loc^+(\R_+^{d+1}\X(-L,L))$, where $\M_\loc^+(\R_+^{d+1}\X(-L,L))$ denotes the space of non-negative Radon measures on $\R_+^{d+1}\X(-L,L)$, with locally finite total variation.  Indeed, \eqref{e1.6}  and \eqref{e1.7} imply, as usual, that for any convex function $\eta:\R\to\R$, we have, in the sense of the distributions in $\R_+^{d+1}$, 
 \begin{equation}\label{e2.2}
 \eta(u)_t+\nabla\cdot \bff_{\eta}(u)-\Delta b_{\eta}(u)=-m_{\eta}(t,x),
 \end{equation}
 where $\bff_\eta'(u)=\eta'(u)\bff'(u)$, $b_\eta'(u)=\eta'(u)b'(u)$, and $m_\eta\in\M_\loc^+(\R^{d+1})$.  Since, we have
 \begin{gather}
 \eta(u)=\int_{\R}\eta'(\xi)\chi(\xi;u)\,d\xi,\quad \bff_{\eta}(u)=\int_{\R}\eta'(\xi)\bff'(\xi)\chi(\xi;u)\,d\xi, \label{e2.3}\\
  b_\eta(u)=\int_{\R}\eta'(\xi)b'(\xi)\chi(\xi;u)\,d\xi, \notag
 \end{gather}
 and, for each fixed $\xi\in\R$,  $\eta_\xi(u)=(\xi-u)_+-(\xi)_+$, with $(\s)_+=\max\{\s,0\}$,  is convex,  writing \eqref{e2.2} for $\eta=\eta_\xi$,  deriving with respect to $\xi$, multiplying by $\eta'(\xi)$, for an arbitrary smooth convex function $\eta$, and integrating the result with respect to $\xi$, we conclude that $m_\eta(t,x)=\int_{\R}\eta''(\xi)\,m(t,x,\xi)$, with $m(t,x,\xi)=m_{\eta_\xi}(t,x)$, and $f(t,x,\xi)=\chi(\xi;u(t,x))$ satisfies \eqref{e2.1},  in $\DD'(\R_+^{d+1}\X(-L,L))$.   

Equation \eqref{e2.2} implies that for any convex entropy $\eta$, the vector field $F=(\eta(u), \bff_\eta(u)-\nabla b_\eta(u))\in \DM^2(C_{R,T})$, for any $R>0$, $T>0$, that is, it is an $L^2$-divergence-measure field on $C_{R,T}$.  By theorems~3.1 and 3.2 in \cite{Fr}, the normal trace of the $\DM^2$-field $F$ at the hyperplane $t=t_*\in(0,T)$, from above, that is, as a part of the boundary of $C_{R,T}\cap \{t>t_*\}$, as well as from below, that is, as part of the boundary of $C_{R,T}\cap\{t<t_*\}$, is simply given by
$$
\la F\cdot\nu, \phi\ra_{t=t_*\pm}=\int_{\R^d}\eta(u(t_*,x))\phi(x)\,dx,
$$
for a.e.\ $t_*>0$, for any $\phi\in C_c^1(\R^d)$, where $\la F\cdot\nu,\cdot\ra_{t=t_*+}$ denotes the normal trace at $\{t=t_*\}$ from above and  $\la F\cdot\nu,\cdot\ra_{t=t_*-}$ the one from below. 
Also, from theorem~3.2 in \cite{Fr}, we deduce that, for any $t_0>0$, 
$$
\la F\cdot\nu, \phi\ra_{t=t_0\pm}=\esslim_{t\to t_0\pm} \int_{\R^d} \eta(u(t,x))\phi(x)\,dx,
$$
for any $\phi\in C_c^1(\R^d)$, and for $t_0=0$ we have, similarly,
$$
\la F\cdot\nu, \phi\ra_{t=0}=\esslim_{t\to 0+} \int_{\R^d} \eta(u(t,x))\phi(x)\,dx.
$$
Now, using \eqref{e2.3} for an arbitrary convex $\eta$,  we deduce that, for $f(t,x,\xi)=\chi(\xi;u(t,x))$,  there exists the limit
\begin{equation}\label{e2.4}
\lim_{t\to t_0+} f(t,\cdot,\cdot)=f^\tau(\cdot,\cdot),
\end{equation}
in the weak star topology of $L^\infty(C_R\X (-L,L))$, for any $R>0$, and any $L>0$ satisfying $\|u\|_{L^\infty(\R_+^{d+1})}\le L$.  Similarly, we have
\begin{equation}\label{e2.4'}
\lim_{t\to t_0-} f(t,\cdot,\cdot)=f^\tau(\cdot,\cdot),
\end{equation}
in the weak star topology of $L^\infty(C_R\X(-L,L))$, if $t_0>0$, since the normal traces from above and below must coincide, as an easy consequence of the Gauss-Green Formula ({\em cf.} \cite{Fr}).  

Following the method in \cite{Va}, in order to prove that the limits in \eqref{e2.4} and \eqref{e2.4'} can be taken as the strong convergence in $L^1(C_{R,T}\X(-L,L))$, it suffices to prove that $f^\tau(\cdot,\cdot)$ is a $\chi$-function, which is proved by using localization method introduced in \cite{Va}. For simplicity we just consider the case $t_0=0$.

Fixing, $x_0\in\R^d$, we consider the sequence
$$
f_\ve(\ut,\ux ,\xi):= f(\ve\ut, x_0+\Lambda(\ve) \ux,  \xi),
$$ 
where $\Lambda(\ve)\ux=(\ve \ux', \ve^{1/2} \ux'')$. So, $f_\ve$ satisfies
\begin{equation} \label{e2.5}
 \po_{\ut} f_\ve+\bfa(\xi)'\cdot\nabla_{\ux'} f_\ve+\ve^{1/2}\bfa(\xi)''\cdot\nabla_{\ux''} f_\ve -b'(\xi)\Delta_{{\ux}''} f_\ve=\po_{\xi} m_\ve,
 \end{equation}
where $\bfa(\xi)'=(\pi_{d'}(\bfa(\xi),\underset{d''}{\underbrace{0,\cdots,0}})$,  $\bfa(\xi)''=\bfa(\xi)-\bfa(\xi)'$, and $m_\ve\in \M_\loc^+((0,\infty)\X\R^d\X\R)$ is defined,  for every $0\le R_1^0< R_2^0$, $R_1^i<R_2^i$, $i=1,\dots,d$, $L_1<L_2$, by  
\begin{multline}
m_\ve\left(\prod_{0\le i\le d} [R_1^i,R_2^i]\X[L_1,L_2]\right)\\
=\frac{1}{\ve^{d'+\frac{d''}2}}m\left([\ve R_1^0, \ve R_2^0]\X(x_0+\Lambda(\ve)\prod_{1\le i\le d} [R_1^i,R_2^i])\X[L_1,L_2]\right),
\end{multline}
where $\Lambda(\ve):\R^d\to\R^d$ is defined by $\Lambda(\ve)z:=(\ve z',\ve^{1/2} z'')$. Following \cite{Va}, as in \cite{FL}, we have there exists a sequence $\ve_n$ converging to 0 and a set $\E\subset \R^d$, with $\LL^d(\R^d\setminus\E)=0$, such that for all $x_0\in\E$ 
\begin{equation}\label{e2.6} 
\lim_{\ve\to0} m_{\ve_n}=0,
\end{equation}
in the weak topology of $\M_{\loc}^+((0,\infty)\X\R^d\X \R)$. 

We now observe that 
\begin{equation}\label{e2.7}
f_\ve(0,\ux,\xi) =f^\tau(x_0+\Lambda(\ve)\ux,\xi).
\end{equation}
Again following \cite{Va}, as in \cite{FL},  we have that there exists a subsequence still denoted $\ve_n$ and a subset $\E'$ of $\R^d$ such that for every $x_0\in\E'$ and for every $R>0$,
\begin{equation}\label{e2.8}
\lim_{\ve_n\to0}\int_{-L}^L\int_{(-R,R)^d}|f^\tau(x_0,\xi)-f^\tau(x_0+\Lambda(\ve_n)\ux,\xi)|\,d\ux\,d\xi=0.
\end{equation}

 Now, we claim that there exists a sequence $\ve_n$ which goes to 0 and a $\chi$-function $f_\infty\in L^\infty(\R_+\X\R^d\X(-L,L))$ such that $f_{\ve_n}$ converges  strongly to $f_\infty$ in $L_\loc^1(\R_+\X\R^d\X(-L,L))$ and
 \begin{equation}\label{e2.9}
 \po_{\ut} f_\infty+ \bfa(\xi)'\cdot\nabla_{\ux'} f_\infty-b(\xi)\Delta_{\ux''} f_\infty=0.
 \end{equation}

 The proof of the claim is very similar to that of proposition~3 of \cite{Va}, and lemma~3.1 in \cite{FL}, and relies on a particular case of the version of averaging lemma in \cite{PS} (see also \cite{TT}).  Here, we need the following  variation of the standard averaging lemma.
 
 \begin{lemma}\label{L:2.1} Let $N,N',N''$ be positive integers with $N=N'+N''$, $f_n(y,\xi)$ be a  bounded sequence in $L^2(\R^N\X\R)\cap L^1(\R^N\X\R)$,  $\bfg_n^i,\bfg^i\in L^2(\R^N\X\R, \R^{N+1})$ be such that $\bfg_n^i\to \bfg^i$ strongly in $L^2(\R^N\X\R,\R^{N+1})$, $i=1,2$, and for $y\in\R^N$ we write $y=(y',y'')$, $y'\in\R^{N'}$, $y''\in\R^{N''}$. Assume
  \begin{equation}\label{e2.10}
 \a(\xi)'\cdot\nabla_{y'}f_n+\a(\xi)''\cdot\nabla_{y''}f_n-\b(\xi)\Delta_{y''}f_n=\po_\xi\nabla_{y,\xi}\cdot\bfg_n^1+\nabla_{y,\xi}\cdot\bfg_n^2,
 \end{equation}
where $\a(\cdot)'\in C^2(\R;\R^{N'})$, $\a(\cdot)''\in C^2(\R;\R^{N''})$, $\b\in C^2(\R)$ satisfy
 \begin{equation}\label{e2.11}
 \LL^1\bigl\{\xi\in\R\,:\,  \a(\xi)\cdot\z'=0\}=0,\qquad \text{for every $\z'\in\R^{N'}$, with $|\z'|=1$},
 \end{equation} 
 where  $\LL^1$ is the Lebesgue measure on $\R$, and also 
 \begin{equation}\label{e2.12}
 \LL^1\{\xi\in\R\,:\, \b(\xi)=0\}=0.
 \end{equation}
Then, for any $\phi\in C_c^\infty(\R)$, the average $u_n^\phi(y)=\int_{\R}\phi(\xi)f_n(y,\xi)\,d\xi$ is relatively compact in $L^2(\R^N)$. 
\end{lemma}

The application of Lemma~\ref{L:2.1} to the problem at hand is made, as in \cite{Va}, by multiplying \eqref{e2.5} by $\phi_1(\ut,\ux),\phi_2(\xi)$ where $\phi_1\in C_0^\infty((1/(2R),2R)\X(-2R,2R)^d)$, $\phi_2\in C_0^\infty(-2L,2L)$, both taking values in $[0,1]$, with $\phi_1(\ut,\ux)=1$, for $(\ut,\ux)\in (1/R,R)\X(-R,R)^d$,  $\phi_2(\xi)=1$, for $\xi\in(-L,L)$. We then consider the equation obtained for $\phi_1\phi_2 f_\ve$, which is easily seen to satisfy the hypotheses of Lemma~\ref{L:2.1}, we refer to \cite{Va} for the details.      

The final step of the proof is to proof that for every $x_0\in\E'$, 
\begin{equation}\label{e2.13}
f_\infty(0,\ux,\xi)=f^\tau(x_0,\xi),
\end{equation}
for a.e.\ $(\ux,\xi)\in\R^d\X(-L,L)$, which the result corresponding to proposition~4 of \cite{Va}. The proof is the same as the one of the referred proposition, and consists in proving that, for any $\phi\in C_0^\infty(\R^d\X(-L,L))$, the sequence 
$$
h^\ve_\phi(t):=\int_{-L}^L\int_{\R^d}(f_\ve(\ut,\ux,\xi)-f_\infty(\ut,\ux,\xi))\phi(\ux,\xi)\,d\ux\,d\xi,
$$
converges to 0 in $BV((0,1))$, which is done exactly as in \cite{Va}.  

Finally, from \eqref{e2.9} and \eqref{e2.13}, we easily conclude that 
$$
f_\infty(\ut, \ux,\xi)=f^\tau(x_0,\xi),
$$
for almost all $(\ut,\ux,\xi)\in \R^{d+1}\X(-L,L)$, which is constant with respect to $(\ut,\ux)$. Hence, since $f_\infty$ is a $\chi$-function for almost all $(\ut,\ux)$, we conclude that 
$f^\tau(x_0,\cdot)$ is a $\chi$-function, as was to be proved. The proof of the strong trace property at any hyperplane $\{t=t_0\}$, $t_0>0$, both from above and from below, follows exactly as just done for $t_0=0$, from above, and this concludes the proof of Theorem~\ref{T:1.1}.

\section{Well-posedness and decay of almost periodic solutions.}\label{S:3}

In this section we are going to prove Theorem~\ref{T:1.2}. We first address the more standard part which is the one asserting the existence and uniqueness of solutions, which we state separately in two lemmas. However, in order to prove  

Since $u_0\in L^\infty(\R^d)$, we may define $u_{\max}=\sup_{\R^d} u_0(x)$ and $u_{\min}=\inf_{\R^d} u(x)$, and we may assume $u_{\max}-u_{\min}>0$, since for $u_0$ constant the problem is trivial.

\begin{lemma}\label{L:3.1}{\bf (Existence)} There exists a weak entropy solution to the problem \eqref{e1.1},\eqref{e1.2}.
\end{lemma}

\noindent
\textsc{Proof:} To prove the existence of weak entropy solution to \eqref{e1.1}, \eqref{e1.2},   we consider the following approximate problems.  For any $R>0$,  let $B_R=B(0,R)$,  be the open ball with radius $R$, centered at the origin, in $\R^d$.  We consider the initial-boundary value problem in $(0,\infty)\X B_R$,    given by
\begin{align}
& u_t^{R,\ve}+\nabla_x\cdot\bff(u^{R,\ve})=\Delta_{x''}b(u^{R,\ve}) +\ve \Delta u^{R,\ve},\quad\text{in $(0,\infty)\X B_R$}, \label{e3.1}\\
&u^{R,\ve}=u_{\min}, \quad\text{in $(0,\infty)\X \po  B_R$}, \label{e3.2}\\
&u^{R,\ve}(0,x)=u_0^{R,\ve}(x),\quad \text{in $B_R$}, \label{e3.3}
\end{align}
where $u_0^{R,\ve}$ is a smooth approximation of $u_0^R:=u_0\lfloor B_R$.

The classical solution of the problem \eqref{e3.1}-\eqref{e3.3} follows from the results in \cite{LSU}.  Multiplying \eqref{e3.1} by $\sgn(u-k) \phi$, for $0\le \phi\in C_0^\infty((-\infty,\infty)\X B_R)$, through standard arguments (e.g., \cite{Kr}), we obtain that the classical solution of \eqref{e3.1}, $u^{R,\ve}$,  satisfies
\begin{multline}\label{e3.4}
\iintl_{(0,\infty)\X  B_R}|u^{R,\ve}-k|\phi_t+\sgn(u^{R,\ve}-k)(\bff(u^{R,\ve})-\bff(k))\cdot\nabla \phi\\
 - (\nabla_{x''}|b(u^{R,\ve})-b(k)|+\ve\nabla|u^{R,\ve}-k|)\cdot\nabla\phi\,dx\,dt \\+\intl_{B_R}|u_0^{R,\ve}(x)-k|\phi(0,x)\,dx\ge0.
\end{multline} 
Clearly also for any $\varphi\in C_0^\infty((-\infty,\infty)\X B_R)$, we have
\begin{multline}\label{e3.5}
\iintl_{(0,\infty)\X  B_R} u^{R,\ve}\varphi_t+ \bff(u^{R,\ve})\cdot\nabla\varphi\,dx\,dt -(\nabla_{x''}b(u^{R,\ve})+\ve\nabla u^{R,\ve})\cdot\nabla\varphi\\+\intl_{B_R}u_0^{R,\ve}(x)\varphi(0,x)\,dx=0.
\end{multline}
By the maximum principle we also deduce that $u_{\min}\le u^{R,\ve}\le u_{\max}$. 

Also, multiplying \eqref{e3.1} by $(b(u^{R,\ve})-b(u_{\min}))$, integrating over $(0,T)\X B_R$, using integration by parts, we deduce
\begin{equation}\label{e3.6}
\iintl_{(0,T)\X B_R} |\nabla_{x''} b(u^{R,\ve})|^2\,dx\,dt \le \intl_{B_R}(|B(u^{R,\ve})|+|B(u_0^{R,\ve})|)\,dx,
\end{equation}
where $B(u)=\int_{u_{\min}}^u (b(z)-b(u_{\min}))\, dz$, and we use the fact that $\bff_b(u)=\int_{u_{\min}}^u (b(z)-b(u_{\min}))\bff'(z)\, dz$ vanishes at the boundary.  It  then follows that
\begin{equation}\label{e3.7}
 \iintl_{(0,T)\X B_R} |\nabla_{x''} b(u^{R,\ve})|^2\,dx\,dt \le C R^d,
 \end{equation}
 for some constant $C>0$ depending only on the data of the problem.  Also, for $R_0+1<R$, multiplying \eqref{e3.1} by $b(u^{R,\ve})\rho_{R_0}$, where $\rho_{R_0}\in C_0^\infty(\R^d)$, $0\le \rho_{R_0}(x)\le 1$, for all $x\in\R^d$, and $\rho_{R_0}(x)=1$, $x\in B_{R_0}$, and $\rho_{R_0}(x)=0$, for $x\in \R^d\setminus B_{R_0+1}$, after standard manipulations, we arrive at
\begin{equation}\label{e3.7'}
 \iintl_{(0,T)\X B_{R_0}} |\nabla_{x''} b(u^{R,\ve})|^2\,dx\,dt \le \tilde C (R_0+1)^d,
 \end{equation}  
 where $\tilde C>0$ only depends on the data of the problem. 
 
 Now, by using the kinetic formulation like \eqref{e2.1} and the averaging lemma as stated in Lemma~\ref{L:2.1}, with a slight modification by including the additional term $\ve\Delta f^\ve$, which does not cause any harm, we easily deduce that, for fixed $R$, when $\ve\to0$, the sequence $u^{R,\ve}$ is relatively compact in $L_{\loc}^1((0,\infty)\X B_R)$, so we may extract a subsequence $u^{R,\ve_n}$ converging in $L_{\loc}^1((0,\infty)\X B_R)$ to a function $u^R$. Clearly, $u^R$ satisfies
\begin{multline}\label{e3.8}
\iintl_{(0,\infty)\X  B_R}|u^{R}-k|\phi_t+\sgn(u^{R}-k)(\bff(u^{R})-\bff(k))\cdot\nabla \phi\\
 - \nabla_{x''}|b(u^{R})-b(k)|\cdot\nabla_{x''}\phi\,dx\,dt \\+\intl_{B_R}|u_0^{R}(x)-k|\phi(0,x)\,dx\ge0,
\end{multline}    
for all $0\le \phi\in C_0^\infty((-\infty,\infty)\X B_R)$, and also
\begin{multline}\label{e3.9}
\iintl_{(0,\infty)\X  B_R} u^{R}\varphi_t+ \bff(u^{R})\cdot\nabla\varphi -\nabla_{x''}b(u^R)\cdot\nabla_{x''}\varphi \,dx\,dt \\+\intl_{B_R}u_0^{R,\ve}(x)\varphi(0,x)\,dx=0,
\end{multline}
for all  $\varphi\in C_0^\infty((-\infty,\infty)\X B_R)$. Moreover, since $\nabla_{x''} b(u^{R,\ve_n})$ weakly converges to $\nabla_{x''}b(u^R)$ in $L_{\loc}^2((0,\infty)\X B_R)$ as is easily seen, from Fatou's lemma we get
\begin{equation}\label{e3.10}
 \iintl_{(0,T)\X B_R} |\nabla_{x''} b(u^{R})|^2\,dx\,dt \le C R^d,
 \end{equation}
for some constant $C>0$ depending only on the data, and from \eqref{e3.7'}
\begin{equation}\label{e3.10'}
 \iintl_{(0,T)\X B_{R_0}} |\nabla_{x''} b(u^{R})|^2\,dx\,dt \le C (R_0+1)^d,
 \end{equation}  

Finally, we send $R\to\infty$, and apply again the kinetic formulation \eqref{e2.1} the averaging lemma as as stated in Lemma~\ref{L:2.1}, to obtain that the sequence $u^R$ is
relatively compact in $L_{\loc}((0,\infty)\X\R^d)$, and so we may extract a subsequence $u^{R_n}$ converging in $L_{\loc}^1((0,\infty)\X\R^d)$ to a function $u$. We then easily deduce that $u\in L^\infty((0,\infty)\X\R^d)$, actually $u_{\min}\le u(t,x)\le u_{\max}$, a.e.\ $(t,x)\in(0,\infty)\X\R^d$, and it satisfies \eqref{e1.6} and \eqref{e1.7}. Also, since $\nabla_{x''}b(u^{R_n})$ weakly converges in $L_{\loc}^2((0,\infty)\X\R^d)$ to $\nabla_{x''}b(u)$,  applying  by Fatou's lemma to \eqref{e3.10'}, we deduce that $u$ satisfies also \eqref{e1.5} and \eqref{e1.5'}, which finishes the proof.

$\hfill\square$

In the next lemma we prove the uniqueness of the weak entropy solution of \eqref{e1.1},\eqref{e1.2}. 

\begin{lemma}\label{L:3.2}{\bf (Uniqueness)} The problem \eqref{e1.1},\eqref{e1.2} has at most one weak entropy solution.
\end{lemma}

\noindent
\textsc{Proof:}  The proof follows through standard arguments ({\em cf.}, e.g., \cite{VH}). So, let $u,v\in L^\infty(\R_+^{d+1})$ be two weak entropy solutions. As in the proof of Proposition~\ref{P:1.1},  by using the doubling of variables method of Kruzhkov \cite{Kr}, as adapted by Carrillo \cite{Ca} to the isotropic degenerate parabolic case, we obtain
\begin{equation}\label{e3.11}
\iintl_{\R_+^{d+1}}|u-v|\phi_t+ \sgn(u-v)(\bff(u)-\bff(v))\cdot\nabla\phi+ |b(u)-b(v)|\Delta_{x''}\phi\,dx\,dt\ge0,
\end{equation}
for all $0\le\phi\in C_0^\infty(\R_+^{d+1})$. We take $\phi(t,x)=\rho(x)\chi_\nu(t)$, where $\rho(x)=e^{-\sqrt{1+x^2}}$ and $\chi_\nu$ is as in the proof of Proposition~\ref{P:1.1}.
We observe that 
$$
\sum_{i=1}^d|\po_{x_i}\rho(x)|+ \sum_{i,j=1}^d|\po_{x_i x_j}^2\rho(x)|\le C\rho(x),
$$
for some constant $C>0$ depending only on $d$. Hence, making $\nu\to0$, we arrive at
\begin{multline*}
\intl_{\R^d}|u(t_1,x)-v(t_1,x)|\rho(x)\,dx\le \int_{\R^d}|u(t_0,x)-v(t_0,x)|\rho(x)\,dx\\+\tilde C \int_{t_0}^{t_1}\int_{\R^d}|u(s,x)-v(s,x)|\rho(x)\,dx\,dt,
\end{multline*}
for a.e. $0<t_0<t_1$, for some $\tilde C>0$ depending only on $\bff,b$ and $d$. Therefore, using Gronwall and the strong trace at $t=0$, we conclude
\begin{equation}\label{e3.12}
\intl_{\R^d}|u(t,x)-v(t,x)|\rho(x)\,dx\le e^{\tilde C t}\int_{\R^d}|u_0(x)-v_0(x)|\rho(x)\,dx,
\end{equation}
which gives the desired result.
 
$\hfill\square$
  
In the next lemma, we prove that the solution operator for \eqref{e1.1},\eqref{e1.2} take  bounded Besicovitch almost periodic functions into bounded Besicovitch almost periodic functions and that $\Gr(u(t,\cdot))\subset\Gr(u_0(\cdot))$.

\begin{lemma}\label{L:3.3}  Let $u(t,x)$ be the weak entropy solution of \eqref{e1.1},\eqref{e1.2}  with $u_0$ satisfying \eqref{e1.u0}. Let $G_0=\Gr(u_0)$. 
Then, $u(t,x)\in C([0,\infty),\BAP(\R^n))\cap L^\infty(\R_+^{d+1})$    and $\Sp(u(t,\cdot))\subset G_0$, for all $t>0$. 
\end{lemma}

\noindent
\textsc{Proof:} The proof follows by the elegant method of reduction to the periodic case introduce by Panov in \cite{Pv2}, more specifically theorems~2.1 and 2.2 in \cite{Pv2}.  Here we limit ourselves to indicate the few adaptations that need to be made. The method begins by considering the case where the initial function $u_0$ is given by a trigonometric polynomial,
\begin{equation}\label{e3.trig}
u_0(x)= \sum_{\l\in\Lambda} a_\l e^{2\pi i\l\cdot x},
\end{equation}
where $\Lambda=\Sp(u_0)\subset\R^d$ is a finite set. Since $u_0$ is real we have that $-\Lambda=\Lambda$  and $a_{-\l}=\bar a_\l$, where as usual $\bar z$ is the complex conjugate of $z\in\C$.   The first observation is that we may find a basis for $G_0$, $\{\l_1,\cdot,\l_m\}$,  so that any $\l\in G_0$ can be uniquely written as $\l=\l(\bar k)=\sum_{j=1}^m k_j\l_j$, $\bar k=(k_1,\cdots,k_m)\in\Z^m$, and the vectors $\l_j$ are linearly independent over $\Z$ and so also over $\Q$. Let $J=\{\bar k\in\Z^m\,:\,\l(\bar k)\in\Lambda\}$.
Then 
\begin{equation}\label{e3.trig'}
u_0(x)=\sum_{\bar k\in J} a_{\bar k} e^{2\pi i\sum_{j=1}^m k_j\l_j\cdot x},\quad a_{\bar k}:=a_{\l(\bar k)}.
\end{equation}
We then have $u_0(x)=v_0(y(x))$, where
\begin{equation}\label{e3.v0}
v_0(y)=\sum_{\bar k\in J}a_{\bar k} e^{2\pi\bar k\cdot y}
\end{equation}
is a periodic function, $v(y+e_i)=v(y)$, $i=1,\cdots,m$, $e_i$ the elements of the canonical basis of $\R^m$, and
$$
y(x)=(y_1,\cdots,y_m),\quad y_j=\l_j\cdot x=\sum_{k=1}^d\l_{jk} x_k,\quad \l_j=(\l_{j1},\cdots,\l_{jd}).
$$  
We then consider the nonlinear degenerate parabolic-hyperbolic equation
\begin{equation}\label{e3.13}
v_t+ \nabla_y\cdot \tilde\bff(v)= (A\nabla_y)\cdot (A\nabla_y) b(v),\quad v(t,y), \ t>0,\ y\in\R^m,
\end{equation}
with $\tilde \bff=(\tilde f_1,\cdots,\tilde f_m)$ and
$$
\tilde f_j(v)=\l_j\cdot\bff (v)=\sum_{k=1}^d \l_{jk} f_k(v),\quad j=1,\cdots,m, \quad A=\frac{\po y}{\po x''}^\top,
$$
and
$$
A\nabla_y=\frac{\po y}{\po x''}^T\nabla_y= (\underset{d'}{\underbrace{0,\cdots,0}}, \sum_{j=1}^m \l_{j(d'+1)}\po_{y_j},\cdots,\sum_{j=1}^m \l_{jd}\po_{y_j}).
$$
We consider the  Cauchy problem for \eqref{e3.13} with initial data
\begin{equation}\label{e3.14}
v(0,y)=v_0(y).
\end{equation}
We say that $v(t,y)\in L^\infty(\R_+^{m+1})$ is an entropy solution of \eqref{e3.13},\eqref{e3.14} if  
\begin{equation}\label{e3.15}
A\nabla_y v\in L_{\loc}^2([0,\infty)\X\R^m),
\end{equation}
and, for all $0\le\tilde\phi\in C_0^\infty(\R^{m+1})$, $k\in\R$, we have
\begin{multline}\label{e3.16}
\iintl_{\R_+^{m+1}} |v-k|\tilde\phi_t+\sgn(v-k)(\tilde\bff(v)-\tilde\bff(k))\cdot\nabla_y\tilde\phi-A\nabla_y|b(v)-b(k)|\cdot A\nabla_y\tilde\phi \,dy\,dt \\
+\intl_{\R^m}|v_0(y)-k|\,\tilde\phi(0,y)\,dy\ge0.
\end{multline}
Existence and uniqueness of  the entropy solution $v(t,y)\in L^\infty(\R_+^{m+1})$ of  \eqref{e3.13},\eqref{e3.14} foliows by the same arguments in the proofs of  Lemmas~\ref{L:3.1} and \ref{L:3.2} and it is easy to see that $v(t,y)$ is also spatially periodic, namely, $v(t,y+e_i)=v(t,y)$, for all $y\in\R^m$, $t>0$, where $e_j$, $j=1,\cdots,m$, is the canonical basis of $\R^m$.  The following assertion corresponds to theorem~2.1 of \cite{Pv2} and it is proof  follows by the same line as the proof of  that result, so we just refer to \cite{Pv2} for the proof.

{\em Assertion \#1.  For a.e.\ $z\in \R^m$ the function $u(t,x)=v(t,z+y(x))$ is a weak  entropy solution of \eqref{e1.1},\eqref{e1.2} with initial data $v_0(z+y(x))$.}   

The next step is another observation in \cite{Pv2} that it follows from Birkhoff individual ergodic theorem \cite{DS} that, for any $w\in L^1(\Pi^m)$, where $\Pi^m:= \R^m/\Z^m$, for almost all $z\in\Pi^m$,  we have
\begin{equation}\label{e3.17}
\Medint_{\R^m}w(z+y(x))\,dx=\int_{\Pi^m} w(y)\,dy.
\end{equation}
Moreover, if $w\in C(\Pi^m)$, then \eqref{e3.17} holds for all $z\in\Pi^m$ and $w^z(x):=w(z+y(x))$ is a (Bohr) almost periodic function for each $z\in\Pi^m$. 

The next main assertion corresponds  to the first part of theorem~2.2 of \cite{Pv2}, that is, it does not include the part about the  decay of the entropy solution, and again its proof follows exactly as the one of the referred theorem and we refer to \cite{Pv2} for the proof.  We leave the claim about the decay of the weak entropy solution to be addressed in a subsequent statement by itself.

{\em Assertion \#2.  Let $u(t,x)$ be a weak entropy solution of \eqref{e1.1},\eqref{e1.2}, assume that the initial function $u_0(x)$ is a trigonometrical polynomial with $G_0=\Gr(u_0)$. Then 
$$
u\in C([0,\infty), \BAP(\R^d))\bigcap L^\infty(\R_+^{d+1})
$$
and $\Sp(u(t,\cdot))\subset G_0$ for all $t>0$.}

We just observe that Assertion~\#2 is proved ({\em cf.} \cite{Pv2}) by using Assertion~\#1 and showing, for a suitable sequence $z_l$ converging to 0 as $l\to\infty$, belonging to the set of full measure of $z\in\R^m$ given by Assertion~\#1, for each fixed $t$ in a set of full measure in $\R_+$, the convergence of the entropy solutions $u^{z_l}(t,\cdot)=v(t,z_l+y(x))$ in $\BAP(\R^d)$, as $z_l\to0$, uniformly with respect to $t$, and  using that for each $z_l$  
$$
u^{z_l}\in C([0,\infty), \BAP(\R^d))\bigcap L^\infty(\R_+^{d+1})
$$
and $\Sp(u^{z_l}(t,\cdot))\subset G_0$ for all $t>0$, which follows from the fact that $v(t,z_l+y)\in C([0,\infty), L^1(\Pi^m))$,  as a consequence of the strong trace property Theorem~\ref{T:1.1}, and the fact that $u^{z_l}(t,x)=v(t,z_l+y(x))$.

Now, let us consider the general case where $u_0\in \BAP(\R^d)\cap L^\infty(\R^d)$. Let $u(t,x)$ be the weak entropy solution of \eqref{e1.1},\eqref{e1.2} obtained above. Following \cite{Pv2},  let $\Gr(u_0)$ be the minimal additive subgroup of $\R^d$ containing $\Sp(u_0)$. We then consider a sequence $u_{0l}$ of trigonometrical polynomials such that $u_{0l}\to u_0$ as $l\to\infty$, in $\BAP(\R^d)$ and $\Sp(u_{0l}\subset \Gr(u_0)$, which may be obtained from the Bochner-Fej\'er trigonometrical polynomials  (see \cite{B}, p.105). We denote by $u_l(t,x)$ the weak entropy solution of \eqref{e1.1},\eqref{e1.2} with initial function $u_{0l}(x)$. By Proposition~\ref{P:1.1}, there exists a set $F\subset \R_+$ of full measure such that, for all $t\in F$ and for every $l\in\N$, we have
\begin{equation}\label{e3.18}
N_1(u(t,\cdot)-u_l(t,\cdot))\le N_1(u_{0l}-u_0)\to 0, \quad\text{as $l\to\infty$}.
\end{equation}
Since $u_{0l}$ has finite spectrum, by Assertion~\#2 we see that  
$$
u_l(t,x)\in C([0,\infty),\BAP(\R^d))
$$ 
and $\Sp(u_l(t,\cdot))\subset \Gr(u_0)$, for all $l\in\N$. Since \eqref{e3.8} holds uniformly with respect to $t\in F$ we conclude that $u(t,\cdot)$ may be extended from $F\cap [0,\infty)$ as a function in $C([0,\infty),\BAP(\R^d))$. Moreover, we easily see that $\Sp(u(t,\cdot))\subset \Gr(u_0)$, for all $t>0$. 

$\hfill\square$

We now turn to the proof of the decay property \eqref{e1.8'}. 

\begin{lemma}\label{L:3.4} The weak entropy solution of \eqref{e1.1},\eqref{e1.2}, with initial function $u_0$ verifying \eqref{e1.u0}, satisfies  \eqref{e1.8'}.
\end{lemma}

\noindent
\textsc{Proof:} Again, reasoning as in \cite{Pv2}, consider first the case where $u_0$ is a trigonometrical polynomial as in \eqref{e3.trig} and \eqref{e3.trig'}, which may be a Bochner-Fej\'er trigonometric polynomial  for the original initial function,  and let $v(t,y)$ be the entropy solution of \eqref{e3.13},\eqref{e3.14}. Writing $A\nabla_y\cdot A\nabla_y b(v)$ in the form 
$\nabla_y B(v)\nabla_y u$, with $A=\frac{\po y}{\po x''}^\top$ as above, we have
$$
B(v)= b'(v) A^\top A,
$$
 so that, for any vector $\k\in\R^m$, we have
 $$
 \k^T B(\xi) \k= b'(\xi) |A\k|^2.
 $$
 Also,
 $$
 \tilde\bff'(v)\cdot \k=\bff'(v)\cdot \tilde A\k, \quad \tilde A=\frac{\po y}{\po x}^\top.
 $$ 
Therefore, we see that conditions \eqref{e1.3} and \eqref{e1.4} imply that, for all $(\tau,\k)\in\R^{m+1}$, with $\tau^2+\k^2=1$, 
\begin{equation}\label{e3.19}
\LL^1\{\xi\in\R\,:\, |\xi|\le \|v_0\|_\infty,\; \tau+\tilde\bff'(\xi)\cdot \k=0,\; \k^\top B(v)\k=0\}=0.
\end{equation}
Hence, we can apply the theorem of Chen and Perthame in \cite{CP2} to obtain that
\begin{equation}\label{e3.20}
\lim_{t\to\infty}\int_{\Pi^m}|v(t,y)-\Me(v_0)|\,dy=0,
\end{equation}
where 
$$
\Me(v_0)= \int_{\Pi^m}v_0(y)\,dy=\Medint_{\R^d}u_0(x)\,dx=\Me(u_0).
$$
Now, we use the fact that
\begin{equation}\label{e3.21}
\Medint_{\R^d}|u(t,x)-C|\,dx=\int_{\Pi^m}|v(t,y)-C|\,dy,
\end{equation}
for any $C\in\R$, which follows by using \eqref{e3.17} and an approximation argument, for which we refer to \cite{Pv2}. Taking $C=\Me(u_0)$ and using \eqref{e3.20}, we conclude that \eqref{e1.8'} holds in the case where $u_0$ is a trigonometric polynomial. 

The general case follows easily from the case just proved by using the Bochner-Fej\'er approximation by trigonometric polynomials, which finishes the proof.

$\hfill\square$

Now, as a final remark, we show that the decay property may also be proved using some ideas in \cite{Fr0}.  First, as a consequence of the fact that $\Gr(u(t,x))\subset\Gr(u_0)$, we have the following result. We recall that the space of Stepanoff almost periodic functions (with exponent $p=1$) in $\R^d$, $\SAP(\R^d)$, is defined as the completion of the trigonometric polynomials with respect to the
norm
$$
\|f\|_S :=\sup_{x\in\R^r}\int_{C_1(x)}|f(y)|\,dy=\sup_{x\in\R^d}\int_{C_1}|f(y+x)|\,dy,
$$
where
$$
C_R(x):=\{y\in\R^d\,:\, |y-x|_\infty:=\max_{i=1,\cdots,d}|y_i-x_i|\le R/2\}.
$$
Another characterization of the Stepanoff almost periodic function (S-a.p., for short)  is obtained introducing the concept of $\ve$-period of a function $f$, that is a number $\tau$ satisfying
\begin{equation}\label{e3.eper}
\|f(\cdot+\tau)-f(\cdot)\|_S\le \tau.
\end{equation}
Let $E_S\{\ve,f\}$ denote the set of such numbers. If the set $E_S\{\ve,f\}$ is relatively dense for all positive values of $\ve$, then the function $f$ is S-a.p.\ (see, e.g., \cite{B}). By the set $E_S\{\ve,f\}$ being relatively dense it is meant that there exists a length $l_\ve$, called {\em inclusion interval},  such that for any $x\in\R^d$, $C_{l_\ve}(x)$ contains an element of $E_S\{\ve,f\}$.  The following lemma is of interest in its own. 

\begin{lemma}\label{L:3.5} If $u_0$ is a trigonometric polynomial, then the weak entropy solution of \eqref{e1.1},\eqref{e1.2}, $u(t,x)$, is S-a.p.\ for all $t>0$, and for any $\ve>0$, there exists $l_\ve>0$ which is an inclusion interval for $u(t,\cdot)$, for all $t>0$. 
\end{lemma}

\noindent
\textsc{Proof:}  Clearly, $u_0$, being a trigonometric polynomial, is S-a.p. The fact that $u(t,x)$ is S-a.p.\ for all $t>0$ follows from \eqref{e3.12}, with $v(t,x)=u(t,x+\tau)$ and $\rho(x-x_0)$ instead of $\rho(x)$, from which we deduce
\begin{multline}\label{e3.22}
\intl_{C_1(x_0)}|u(t,x+\tau)-u(t,x)| \,dx\le c(t)\intl_{C_R(x_0)}|u_0(x)-u_0(x+\tau)|\rho(x-x_0)\,dx \\+ O(\frac1{R})
\le c(R,t)\sup_{x\in\R^d}\intl_{C_1(x)}|u_0(y+\tau)-u_0(y)|\,dy+O(\frac1{R}),
\end{multline}
where $c(R,t)$ is a positive constant depending only on $R,t$ and $O(\frac1{R})$ goes to zero when $R\to\infty$ uniformly with respect to $x_0$.  So, choosing $R$ large enough so that $O(1/R)\le \ve/2$ and then taking any $\tau\in E_S\{\ve/(2c(R,t)), u_0\}$, we get that $\tau \in E_S\{\ve, u(t,\cdot)\}$, and so $u(t,\cdot)$ is S-a.p.

Now we prove that for any $\ve>0$, there exists $l_\ve>0$ which is an inclusion interval for $u(t,\cdot)$, for all $t>0$. For this we use the following two results of \cite{B}, p.\ 53, which hold for the classical almost periodic functions (u.a.p., for short) in the line. The $\ve$-periods of  a u.a.p.\ function $f$ are defined by \eqref{e3.eper} with $\|\cdot\|_S$ replaced by the $\sup$-norm $\|\cdot\|_\infty$ and we denote the set of such $\ve$-periods by $E\{\ve,f\}$.

{\em Assertion \#1.  Given a u.a.p.\ function 
\begin{equation}\label{e3.uap}
f(x)\sim \sum_{n=1}^\infty A_n e^{i\Lambda_n x},
\end{equation}
to any positive integer $N$ and a positive number $\d<\pi$ corresponds a positive $\ve$ such that all numbers $\tau$ of the set $E\{\ve, f(x)\}$ satisfy the following Diophantine inequalities:
\begin{equation}\label{e3.DE}
|\Lambda_n\tau | < \d \mod 2\pi,\quad n=1,\cdots,N,
\end{equation}
where the inequality means that there exists $k\in\Z$ such that $|\Lambda_n\tau-2\pi k|<\d$. }

\medskip

{\em Assertion \#2.  Given a u.a.p.\ function as in \eqref{e3.uap}, to any $\ve>0$ corresponds a positive integer $N$ and a positive $\d<\pi$ such that any number $\tau$ satisfying the $N$ Diophantine inequalities \eqref{e3.DE}, belongs to $E\{\ve, f(x)\}$. }

\medskip

These assertions are stated for u.a.p.\ functions in $\R$, but it is easy to see that they can be easily extended to u.a.p.\ functions in $\R^d$.  Now, let us show how these two assertions from \cite{B} can be applied to prove the lemma.   For any $t>0$, we approximate $u(t,x)$ by the corresponding Bochner-Fej\'er trigonometric polynomials, whose spectrum is contained in $\Sp(u(t,\cdot))$ and the coefficients have absolute values dominated by the corresponding coefficient in $u(t,\cdot)$, which, in turn, are dominated by those of the  initial function by Proposition~\ref{P:1.1}. Therefore, Assertion~\#2, whose prove depends only on the frequencies and absolute value of the coefficients ({\em cf.} \cite{B}),  implies that the set of $\tau$'s satisfying \eqref{e3.DE} is contained in $E\{\ve,u(t,\cdot)\}$, for all $t>0$. On the other hand, Assertion~\#1 implies that the set of such $\tau$'s satisfying \eqref{e3.DE} is relatively dense, since it contains $E\{\ve',u_0\}$ for some $\ve'>0$, and such sets are relatively dense. Therefore, we can find $l_\ve>0$ which is an inclusion interval for $u(t,\cdot)$ for all $t>0$, which was to be proved.  

$\hfill\square$

Now we can use Lemma~\ref{L:3.5} to give an alternative proof of the decay property \eqref{e1.8'}. Clearly, from Proposition~\ref{P:1.1}, by approximating the initial function by Bochner-Fej\'er trigonometric polynomials, it suffices to consider the case where the initial function $u_0$ is itself a trigonometric polynomial. Let us then consider the scaling sequence $u^T(t,x):=u(Tt, Tx',\sqrt{T}x'')$, and define $\xi'=x'/t$, $\xi''=x''/\sqrt{t}$.  So, $u^T$ is a uniformly bounded sequence of weak entropy solutions of \eqref{e1.1},\eqref{e1.2}, with initial functions $u_0^T(x):=u_0(Tx',\sqrt{T}x'')$. Using the Averaging Lemma~\ref{L:2.1}, we deduce that $u^T$ is relatively compact in $L_\loc^1(\R_+^{d+1})$ and the initial functions clearly converge weakly to $\bar u_0=\Me(u_0)$. By passing to a subsequence that we still denote by $u^T(t,x)$, we have that $u^T\to \bar u$ as $T\to\infty$, in 
$L^1_\loc(\R^{d+1})$, for some $u\in L^\infty(\R_+^{d+1})$. We see also that $\bar u$ satisfies \eqref{e1.5}, \eqref{e1.5'}, \eqref{e1.6} and \eqref{e1.7}, all of which are easy to be verified, and we observe by \eqref{e1.7} that $\bar u(0,x)=\Me(u_0)$. Now, by the uniqueness Lemma~\ref{L:3.2}, we conclude that $\bar u(t,x)=\Me(u_0)$, that is $u^T\to \Me(u_0)$, in $L_\loc^1(\R^{d+1})$. This, in particular, implies
\begin{multline*}
0=\lim_{T\to\infty}\int_0^1\intl_{|x'|\le c',\, |x''|\le c''}|u(Tt,Tx',\sqrt{T}x'')-\Me(u_0)|\,dx'\,dx''\,dt \\
= \lim_{T\to\infty} \frac1T\int_0^T \frac1{T^{d'+d''/2}}\intl_{|x'|\le c'T,\, |x''|\le c''\sqrt{T}} |u(t,x',x'')-\Me(u_0)|\,dx'\,dx''\,dt\\
\ge \frac1{2^{d'+d''/2}} \lim_{T\to\infty}\frac1T\int_{T/2}^T\intl_{|\xi'|\le c',\ |\xi''|\le c''}|u(t,\xi't,\xi''\sqrt{t})-\Me(u_0)|\,d\xi'\,d\xi''\,dt,
\end{multline*}      
which implies
\begin{equation}\label{e3.23}
\lim_{T\to\infty}\frac1T\int_0^T\intl_{|\xi'|\le c',\ |\xi''|\le c''}|u(t,\xi't,\xi''\sqrt{t})-\Me(u_0)|\,d\xi'\,d\xi''\,dt=0,
\end{equation}
as is easily seen. Now, invoking Lemma~\ref{L:3.5}, we can then make a computation similar to that in p.51 of \cite{Fr0} in order to get, for all $t>0$ large enough, 
\begin{equation}\label{e3.24}
\intl_{|\xi'|\le c',|\xi''|\le c''}|u(t,\xi't,\xi''\sqrt{t})-\Me(u_0)|\,d\xi'\,d\xi''\ge c_1\Me(|u(t,\cdot)-\Me(u_0)|),
\end{equation}
for certain positive constant $c_1$ depending only on the dimension. Therefore, by \eqref{e3.23}, we deduce
\begin{equation}\label{e3.25}
\lim_{T\to\infty} \frac1T \int_0^T \Me(|u(t,\cdot)-\Me(u_0)|)\,dt=0.
\end{equation}
Now, by Proposition~\ref{P:1.1}, we conclude 
 \begin{equation}\label{e3.26}
\lim_{t\to\infty}  \Me(|u(t,\cdot)-\Me(u_0)|)\,dt=0,
\end{equation}
which is the desired result.

\section{Acknowledgements} 

The author gratefully acknowledges the support from CNPq, through grant proc.\ 303950/2009-9, and FAPERJ, through grant proc.\ E-26/103.019/2011.

\end{document}